\documentclass{article}

\usepackage{bm}
\usepackage{amsmath}
\usepackage{amsfonts}
\usepackage{graphicx}

\newcommand{\bsigma}{\bm{\sigma}}
\newcommand{\btau}{\bm{\tau}}
\newcommand{\bn}{\bm{n}}

\newcommand{\bet}{\noalign{\vskip6pt plus 3pt minus 1pt}}
\newcommand{\betmore}{\noalign{\vskip9pt plus 3pt minus 1pt}}

\newcommand{\nn}{\nonumber}

\title{The Compact Discontinuous Galerkin (CDG) Method for Elliptic Problems}

\author{J. Peraire\thanks{Department of Aeronautics and Astronautics,
MIT, 77 Massachusetts Avenue 37-451, Cambridge, MA 02139
(peraire@mit.edu).}
\and
P.-O. Persson\thanks{Department of Mathematics,
MIT, 77 Massachusetts Avenue 2-363A, Cambridge, MA 02139
(persson@mit.edu).}}

\begin{document}

\maketitle

\begin{abstract}
We present a compact discontinuous Galerkin (CDG) method for an
elliptic model problem. The problem is first cast as a system of
first order equations by introducing the gradient of the primal
unknown, or flux, as an additional variable. A standard
discontinuous Galerkin (DG) method is then applied to the
resulting system of equations. The numerical interelement fluxes
are such that the equations for the additional variable can be
eliminated at the element level, thus resulting in a global system
that involves only the original unknown variable. The proposed
method is closely related to the local discontinuous Galerkin
(LDG) method [B.~Cockburn and C.-W.\ Shu, {\it SIAM J.\ Numer.\
Anal.}, 35 (1998), pp.~2440--2463], but, unlike the LDG method,
the sparsity pattern of the CDG method involves only nearest
neighbors. Also, unlike the LDG method, the CDG method works
without stabilization for an arbitrary orientation of the element
interfaces. The computation of the numerical interface fluxes for
the CDG method is slightly more involved than for the LDG method,
but this additional complication is clearly offset by increased
compactness and flexibility. Compared to the BR2 [F.~Bassi and
S.~Rebay, {\it J.\ Comput.\ Phys.}, 131 (1997), pp.~267--279] and
IP [J.~Douglas, Jr., and T.~Dupont, in {\it Computing Methods in
Applied Sciences\/} (Second Internat.\ Sympos., Versailles, 1975),
Lecture Notes in Phys.\ 58, Springer, Berlin, 1976, pp.~207--216]
methods, which are known to be compact, the present method
produces fewer nonzero elements in the matrix and is
computationally more efficient.
\end{abstract}

\section{Introduction}\label{sec1}

Discontinuous Galerkin (DG) methods \cite{RKDG} have become the
subject of considerable research over recent years due to their
potential to overcome some of the perceived shortcomings of the
more established discretization methods. For convection problems,
DG methods produce stable discretizations without the need for
cumbersome stabilization strategies. They work well on arbitrary
meshes and allow for different orders of approximation to be used
on different elements in a very straightforward manner. Clearly,
this flexibility comes at the expense of duplicating the degrees
of freedom at the element boundary interfaces. This is a serious
drawback when low order polynomial approximations are used, but it
is less important for high order interpolations. DG methods appear
to be ideally suited for applications involving wave propagation
phenomena, where low dispersion and high accuracy are required,
such as aeroacoustics or electromagnetics.

While DG methods seem to be well suited for the discretization of
first order hyperbolic problems, their extension to elliptic
problems is far less obvious. A number of extensions to deal with
the elliptic problem have been proposed and analyzed under a
unified framework in~\cite{arnold02unified}. Also, a comparison of
the performance of various schemes from a practical perspective is
presented in~\cite{castillo:524}. Among the various 
alternatives, the local discontinuous Galerkin (LDG) method
\cite{cockburn98ldg} has emerged as one of the most \hbox{popular}
choices. The LDG method appears to be one of the most accurate and
stable schemes among those tested. In addition, the LDG method is
easy to implement for complex convective-diffusive systems and can
be generalized to handle equations involving higher order
derivatives \cite{DGHO}. In the LDG method, the original equation
involving second order derivatives is cast as a system of first
order equations by introducing additional variables for the
solution gradient, or flux. The resulting system is then
discretized using a standard DG approach. By appropriately
choosing the interelement fluxes, the additional variable can be
eliminated locally. Thus, a stable discretization that involves
only the original unknown variable is obtained. Unfortunately,
when the LDG method is used in multiple dimensions, the
discretization generated has the undesirable feature that the
degrees of freedom in one element are connected, not only to those
in the neighboring elements, but also to those in some elements
neighboring the immediate neighbors. For applications employing
explicit or iterative solution techniques, this is usually not a
problem, but for applications where the matrix needs to be formed,
this represents a severe disadvantage.

Two alternative formulations for the treatment of the second order
derivatives are the symmetric interior penalty (IP)
method~\cite{IPMethod} and the BR2 method proposed
in~\cite{bassi97br2}. In these methods, the original form of the
equation involving second derivatives is discretized directly, and
stabilization is added explicitly in a sufficient amount to render
the method stable. Although somewhat simpler, the IP method
appears to be less popular than the BR2 method. This is probably
because of the requirement of a penalty parameter that depends on
both the mesh and the approximation order. Both these methods have
the advantage that they are compact in the sense that only the
degrees of freedom belonging to neighboring elements are connected
in the discretization. When suitable penalization is employed
these approaches are competitive with the LDG scheme in terms of
accuracy. Thus, these schemes are an attractive alternative to the
LDG scheme when an implicit solution of the discretized system is
required.

For many applications of interest involving convective-diffusive
systems, such as the Navier--Stokes equations at high Reynolds
numbers, the time and length scales are such that implicit
discretization turns out to be a requirement. In this paper, we
develop a variation of the LDG method, the compact discontinuous
Galerkin (CDG) method. The main motivation for developing this new
scheme is to eliminate the distant connections between
nonneighboring elements which arise when the LDG scheme is used in
multiple dimensions. We note that in the one-dimensional case the
CDG and LDG schemes are identical, but in the multidimensional
case they differ in the approximation to the solution gradient at
the interface between neighboring elements. This seemingly minor
difference results in a scheme that appears to inherit all the
attractive features of the LDG method and is compact. In addition,
numerical experiments indicate that the CDG scheme is slightly
more stable than the LDG method and is less sensitive to the
element and/or interface orientation. In particular, when the
stabilization constant is set to zero, the CDG scheme is stable in
situations where the LDG method is unstable. It is well known
that, without explicit stabilization, the LDG scheme is stable
only when the orientation of element interfaces satisfies a
certain condition~\cite{sherwin06elliptic}.

Since the CDG scheme is compact, it produces a sparser
connectivity matrix than the LDG scheme, meaning lower storage
requirements and higher computational performance. Thus, the
slight additional increase in complexity involved in the numerical
flux evaluation is more than offset by the increased efficiency
benefits. Compared to the IP and BR2 methods, the CDG scheme is
computationally simpler, generates a sparser matrix with a smaller
number of nonzero elements when using a nodal basis, and appears
to produce slightly more accurate results than the BR2 method in
the numerical tests performed. Given the similarities between the
BR2 and IP methods, we have considered only the BR2 method in our
numerical comparisons.

The remainder of the paper is organized as follows. In
section~\ref{DGFORM}, we introduce our model second order elliptic
problem. Next, we describe the LDG discretization method and adopt
the framework introduced in \cite{arnold02unified} to write the
LDG algorithm in the so-called primal form. This form, involving
only the original problem variable, highlights the symmetry of the
scheme as well as the sparsity pattern. In section~\ref{CDG}, we
present the CDG method. The CDG method is then written in primal
form so that it can be easily compared with the LDG method. Like
the LDG method, the CDG method is shown to be symmetric,
conservative, and adjoint consistent. It turns out that the CDG
and LDG schemes are so closely related that the error estimate
presented in~\cite{arnold02unified} for the LDG method is
essentially applicable to the CDG method without changes. In
section~\ref{STAB}, we compare the LDG and CDG schemes using the
test problem presented in~\cite{sherwin06elliptic}. The increased
stability of the CDG scheme, for arbitrary interface ordering, is
shown numerically by calculating the size of the null-space for
the model test problem. Practical implementation and efficiency
issues such as sparsity patterns and storage requirements for the
LDG, BR2, and CDG schemes, in the more general $d$-dimensional
setting, are addressed in section~\ref{IMPLEMENT}. Finally, we
conclude in section~\ref{NUMERICAL} with some numerical results
aimed at comparing the accuracy and conditioning of the LDG, BR2,
and CDG schemes.

\section{Discontinuous Galerkin formulation}\label{DGFORM}

\subsection{Problem definition}\label{sec2.1}

The proposed method will be described for the model Poisson
problem
\begin{equation}
\begin{array}{rclcll}
- \nabla \cdot ( \kappa \nabla u) & = & f & \qquad & \mbox{in} &
\Omega, \\
\bet
u & = & g_D & \qquad & \mbox{on} & \partial \Omega_D , \\
\bet
\displaystyle \kappa \frac{\partial u}{\partial n} & = & g_N &
\qquad & \mbox{on} & \partial \Omega_N , \end{array}
\label{probdef}
\end{equation}
where $\Omega$ is a bounded domain in $\mathbb{R}^d$ with boundary
$\partial \Omega = \partial \Omega_D \cup \partial \Omega_N$ and
$d = 1,2$, or $3$ is the dimension. Here, $f(\bm x)$ is a given
function in $L^2(\Omega)$, and $\kappa(\bm x) \in
L^\infty(\Omega)$ is positive. Further, we assume that the length
of $\partial \Omega_D$ is not zero.

\subsection{DG Formulation for elliptic problems}\label{sec2.2}

In order to develop a DG method, we rewrite the above
problem (\ref{probdef}) as a first order system of equations
\begin{equation}
\begin{array}{rclcll}
- \nabla \cdot {\bm \sigma} & = & f & \qquad & \mbox{in} & \Omega
, \\
\bet
{\bm \sigma} & = & \kappa \nabla u & \qquad & \mbox{in} & \Omega
, \\
\bet
u & = & g_D & \qquad & \mbox{on} & \partial \Omega_D , \\
\bet
\displaystyle {\bm \sigma} \cdot {\bm n} & = & g_N & \qquad &
\mbox{on} & \partial \Omega_N ,
\end{array} \label{probaux}
\end{equation}
where ${\bm n}$ is the outward unit normal to the boundary of
$\Omega$.

Next, we introduce the {\em broken} spaces $V({\cal T}_h)$ and
$\Sigma({\cal T}_h)$ associated with the triangulation ${\cal T}_h
= \{ K \}$ of $\Omega$. In particular, $V({\cal T}_h)$ and
$\Sigma({\cal T}_h)$ denote the spaces of functions whose
restriction to each element $K$ belongs to the Sobolev spaces
$H^1(K)$ and $[H^1(K)]^d$. That~is,
\begin{eqnarray}
V & = & \{ v \in L^2(\Omega) \ | \ v|_K \in H^1(K) \ \ \forall K \in
{\cal T}_h \}, \\
\bet \Sigma & = & \{ {\btau} \in [L^2(\Omega)]^d \ | \ \btau|_K
\in [H^1(K)]^d\ \ \forall K \in {\cal T}_h \}.
\end{eqnarray}
In addition, we introduce the finite element subspaces $V_h
\subset V$ and $\Sigma_h \subset \Sigma$ as
\begin{eqnarray}
V_h & = & \{ v \in L^2(\Omega) \ | \ v|_K \in {\cal P}_p(K) \ \
\forall K \in {\cal T}_h \}, \\
\bet \Sigma_h & = & \{ {\btau} \in [L^2(\Omega)]^d \ | \ \btau|_K
\in [{\cal P}_p(K)]^d \ \ \forall K \in {\cal T}_h \},
\end{eqnarray}
where ${\cal P}_p(K)$ is the space of polynomial functions of
degree at most $p \ge 1$ on~$K$.

Following~\cite{cockburn98ldg}, we consider DG formulations of the
form: find $u_h \in V_h$ and $\bsigma_h \in \Sigma_h$ such that
for all $K \in {\cal T}_h$ we have
\begin{eqnarray}
\int_K \bsigma_h \cdot \btau \, dx = - \int_K u_h \nabla \cdot
(\kappa \btau) \, dx + \int_{\partial K} \hat{u}
\kappa \btau \cdot \bn \, ds & \quad \forall \btau \in [{\cal P}_p(K)]^d ,& \label{elementDG}\\
\bet
\int_K \bsigma_h \cdot \nabla v \, dx = - \int_K f v \, dx +
\int_{\partial K} \hat{\bsigma} \cdot \bn v \, ds & \forall v \in
{\cal P}_p(K) . & \label{elementDG1}
\end{eqnarray}
Here, the numerical fluxes $\hat{\bsigma}$ and $\hat{u}$ are
approximations to $\bsigma = \kappa \nabla u$ and to $u$,
respectively, on the boundary of the element $K$. The DG
formulation is complete once we specify the numerical fluxes
$\hat{\bsigma}$ and $\hat{u}$ in terms of $\bsigma_h$ and $u_h$
and the boundary conditions.

Expressions (\ref{elementDG}) and (\ref{elementDG1}) apply to each
element separately. In order to write expressions which are
applicable over the whole domain, we require some additional
notation. Here, we closely follow the notation used
in~\cite{arnold02unified}.

Consider two adjacent elements $K^+$ and $K^-$ of the
triangulation ${\cal T}_h$, and denote by $e = \partial K^+ \cap
\partial K^-$ their common face. Further, assume that $\bn^\pm$
denote the unit normals to $\partial K^\pm$, respectively, at any
point on the face $e$. Similarly, let $(\btau^\pm, v^\pm)$ denote
the traces on $e$ of functions $(\btau, v) \in \Sigma_h \times V_h
$ which are smooth in the interior of elements $K^\pm$. The
average and jump operators are given as
\begin{equation}
\nonumber
\begin{array}{lcl}
\{\btau \} = (\btau^+ +
\btau^-)/2, & \hspace{10ex} & \{v\} = (v^+ + v^-)/2, \\
\bet [ \btau ] = \btau^+ \cdot \bn^+ + \btau^- \cdot \bn^-, & &
[v] = v^+ \bn^+ + v^- \bn^-.
\end{array}
\end{equation}
Note that, according to this definition, the jump of a scalar
quantity is a vector, but the jump of a vector quantity becomes a
scalar.

Now, by summing (\ref{elementDG}) and (\ref{elementDG1}) over all
elements and considering only conservative schemes for which the
numerical fluxes $\hat{u}$ and $\hat{\bsigma}$ on a given face are
unique, we obtain the following global expressions: find $u_h \in
V_h$ and $\bsigma_h \in \Sigma_h$ such that
\begin{align}
\int_\Omega \bsigma_h \cdot \btau \, dx & = - \int_\Omega u_h
\nabla_h \cdot (\kappa \btau) \, dx + \int_{{\cal E}_i} \hat{u}
[\kappa \btau] \, ds + \int_{\partial \Omega}
\hat{u}\, \kappa \btau \cdot \bn ds \quad \forall \btau \in \Sigma_h, \label{globalDG}\\
\bet
 \int_\Omega \bsigma_h \cdot \nabla_h v \, dx & = \int_\Omega f v \, dx
 +
\int_{{\cal E}_i} \hat{\bsigma} \cdot [v] \, ds + \int_{\partial
\Omega} v\, \hat{\bsigma} \cdot \bn \, ds  \quad \forall v \in
V_h,
 \label{globalDG1}
 \end{align}
where ${\cal E}_i$ denotes the union of all the interior faces in
the triangulation ${\cal T}_h$. Also, $\nabla_h$ denotes the
broken gradient operator. That is, $\nabla_h v$ and $\nabla_h
\cdot \btau$ are functions whose restriction to $K$ is equal to
$\nabla v$ and $\nabla \cdot \btau$, respectively.

For later use, we note that, if we use the integration by parts
formula,
\begin{equation}
- \int_\Omega v \nabla_h \cdot \btau \, dx = \int_\Omega \btau
\cdot \nabla_h v \, dx - \int_{{\cal E}_i} ([ v] \cdot \{ \btau\}
+ \{ v\} [ \btau] ) \, ds - \int_{\partial \Omega} v \btau \cdot
\bn \, ds,
\end{equation}
which is valid for all $\btau \in [H^1({\cal T}_h)]^d$ and $v \in
H^1({\cal T}_h)$, we can write (\ref{globalDG}) as
\begin{equation}
\begin{split}
\int_\Omega \bsigma_h \cdot \btau \, dx & = \int_\Omega \btau
\cdot (\kappa \nabla_h u_h) \, dx - \int_{{\cal E}_i} ( [u_h] \cdot
\{\kappa \btau\} - \{\hat{u}- u_h \} [\kappa \btau]) \, ds \\
\bet
 & \quad + \int_{\partial \Omega} (\hat{u} - u_h) \kappa \btau \cdot \bn ds \
 \qquad \forall \btau \in \Sigma_h . \label{globalDGibp}
\end{split}
\end{equation}

\subsection{The LDG method}\label{sec2.3}

Since our method is closely related to the LDG method presented in
\cite{cockburn98ldg}, we start with a description of the LDG
algorithm. For the LDG method, the numerical interelement fluxes
$(\hat{\bsigma}, \hat{u})$ are given~by
\begin{eqnarray}
\hat{\bsigma} & = & \{ \bsigma_h \} - C_{11}[u_h] + {\bm C}_{12}
[\bsigma_h], \label{qflux}
\\
\bet \hat{u} & = & \{ u_h \} - {\bm C}_{12} \cdot [u_h]
\label{uflux}
\end{eqnarray}
for the interior faces, and
\begin{equation}
\begin{array}{r@{\,}c@{\,}l@{\,}r@{\,}cll}
 \hat{\bsigma} & = & \bsigma_h - C_{11}(u_h - g_D)\bn, \qquad
\hat{u} & = & g_D
\qquad & \mbox{on} \quad \partial \Omega_D, \\
\bet
 \hat{\bsigma} & = & g_N \bn, \hfill \hat{u} & = & u_h \qquad & \mbox{on}
\quad
\partial \Omega_N,
\end{array} \label{bfluxes}
\end{equation}
for the boundary faces. Here, $C_{11}$ is a positive constant and
${\bm C}_{12}$ is a vector which is determined for each interior
face according to
\begin{equation}
 {\bm C}_{12} = \frac{1}{2} ( S_{K^+}^{K^-} \bn^+ + S_{K^-}^{K^+}
 {\bn^-}), \label{switch}
\end{equation}
where $S_{K^+}^{K^-} \in \{ 0, 1 \}$ is a {\em switch} which is
defined for each element face. That is, $S_{K^+}^{K^-}$ denotes
the switch associated with element $K^+$ on the face that element
$K^+$ shares with element $K^-$. The switches always satisfy that
\begin{equation}
S_{K^+}^{K^-}+ S_{K^-}^{K^+} = 1  \label{sumtoone}
\end{equation}
but are otherwise arbitrary. We note that, although the form
(\ref{switch}) is not the most general form for ${\bm C}_{12}$
presented in~\cite{cockburn98ldg}, other choices lead to wider
stencils in the final discrete equations. We also point out that
the choice of element face switches has an effect on the final
form of the discrete equations.

\subsubsection{Primal form of the LDG algorithm}
\label{sec2.3.1}

In order to derive the primal form of the LDG algorithm, we first
particularize (\ref{globalDGibp}) for the fluxes given by
(\ref{uflux}),
\begin{equation}
\begin{split}
\int_\Omega \bsigma_h \cdot \btau \, dx & = \int_\Omega \btau
\cdot (\kappa \nabla_h u_h) \, dx - \int_{{\cal E}_i} ( [u_h] \cdot
\{\kappa \btau\} + {\bm C}_{12} \cdot [u_h] [\kappa \btau]) \, ds \\
\bet & \quad + \int_{\partial \Omega_D} (g_D - u_h) \kappa \btau
\cdot \bn ds \quad \forall \btau \in \Sigma_h.
\label{globalDGibpLDG}
\end{split}
\end{equation}
To obtain an expression for $\bsigma_h$ as a function $u_h$, we
follow~\cite{arnold02unified} and introduce the lifting operators
$r: [L^2({\cal E}_i)]^d \to \Sigma_h$, $l: L^2({\cal E}_i) \to
\Sigma_h$, and $r_D: L^2(\partial \Omega_D) \to \Sigma_h$:
\begin{equation}
\begin{array}{rcll}
\displaystyle \int_\Omega r(\phi) \cdot \btau \, dx & = &
\displaystyle - \int_{{\cal E}_i} \phi \cdot \{ \btau \} \,
ds  \qquad & \forall \btau \in \Sigma_h,\\
\betmore\betmore \displaystyle \int_\Omega l(q) \cdot \btau \, dx
& = & - \displaystyle \int_{{\cal E}_i} q [ \btau] \, ds  &
\forall
\btau \in \Sigma_h, \\
\betmore\betmore
\displaystyle \int_\Omega r_D(q) \cdot \btau \, dx & = &
 \displaystyle - \int_{\partial \Omega_D} q \btau \cdot \bn \, ds
\qquad & \forall \btau \in \Sigma_h.
\end{array} \label{liftings}
\end{equation}
Thus, we can write (\ref{globalDGibpLDG})~as
\begin{equation}
\begin{gathered}
\int_\Omega (\bsigma_h - \kappa \nabla_h u_h - \kappa r([u_h]) -
\kappa l({\bm C}_{12}\cdot[u_h]) + \kappa r_D(g_D - u_h)) \cdot
\btau \, dx = 0 \\
\bet \forall \btau \in \Sigma_h.
\end{gathered}
\end{equation}
Therefore, we have
\begin{equation}
\bsigma_h = \kappa \nabla_h u_h + \bar{\bsigma}_h,
\label{bsigmah}
\end{equation}
where $\bar{\bsigma} \in \Sigma_h$ is
\begin{equation}
\bar{\bsigma}_h = \kappa r([u_h]) + \kappa l({\bm
C}_{12}\cdot[u_h]) - \kappa r_D(g_D - u_h).
\label{barsigma}
\end{equation}
Thus, we see that that $\bsigma_h$ is equal to $\kappa \nabla_h
u_h$ plus an additional perturbation term which is forced by
$[u_h]$, ${\bm C}_{12} \cdot [u_h]$, and $g_D-u_h$. Also, note
that $r_D(g_D-u_h)$ is nonzero only on the elements that have a
face on the Dirichlet boundary. In writing expressions
(\ref{bsigmah}) and (\ref{barsigma}), we have assumed that
$\nabla_h V_h \subset \Sigma_h$, which is certainly the case if
equal order polynomial interpolants are used for $V_h$
and~$\Sigma_h$.

Setting $\tau = \nabla_h v$ in (\ref{globalDGibpLDG}), we can
rewrite (\ref{globalDG1})~as
\begin{multline}
 \int_\Omega
\nabla_h v \cdot (\kappa \nabla_h u_h) \, dx - \int_{{\cal E}_i} (
[u_h] \cdot \{\kappa \nabla_h v\} + {\bm C}_{12} \cdot [u_h] [\kappa
\nabla_h v]) \, ds \\
\bet
+ \int_{\partial \Omega_D} (g_D - u_h) \kappa
\nabla_h v \cdot \bn ds \\
\bet
= \int_\Omega f v \, dx
 +
\int_{{\cal E}_i} \hat{\bsigma} \cdot [v] \, ds + \int_{\partial
\Omega} v\, \hat{\bsigma} \cdot \bn \, ds  \qquad \forall v \in
V_h. \label{goveeq}
\end{multline}
\noindent Making use of (\ref{qflux}), (\ref{bfluxes}),
(\ref{bsigmah}), and (\ref{barsigma}), the terms involving
$\hat{\bsigma}$ in the above equation can be written~as
\begin{equation}
\begin{split}
 \int_{{\cal E}_i} \hat{\bsigma} \cdot [v] \, ds & = \int_{{\cal E}_i}
( \{\kappa \nabla_h u_h \} + {\bm C}_{12} [ \kappa \nabla_h u_h ])
\cdot [v] \, ds + \int_{{\cal E}_i} (\{ \bar{\bsigma} \} + {\bm
C}_{12} [ \bar{\bsigma} ]) \cdot [v] \, ds \\
\bet
& \quad - \int_{{\cal E}_i} C_{11} [u]\cdot [v] \, ds \\
\bet
 & = \int_{{\cal E}_i} ( \{\kappa \nabla_h u_h \} + {\bm C}_{12} [
\kappa \nabla_h u_h ]) \cdot [v] \, ds \\
\bet
& \quad - \int_{\Omega}
\kappa (r([v]) + l({\bm C}_{12} \cdot [v]))\cdot (r([u_h]) + l({\bm
C}_{12} \cdot [u_h]) + r_D(u_h) ) \, dx \\
\bet
& \quad + \int_{\Omega}
\kappa (r([v]) + l({\bm C}_{12} \cdot [v]))\cdot r_D(g_D) \, dx
 - \int_{{\cal E}_i } C_{11} [u ] \cdot [v] \, ds
 \nonumber
 \end{split}
\end{equation}
and
\begin{equation}
\begin{split}
 \int_{\partial
\Omega} v\, \hat{\bsigma} \cdot \bn \, ds & = \int_{\partial
\Omega_D} v \bsigma_h \cdot \bn \, ds - \int_{\partial \Omega_D}
C_{11} v u_h \, ds + \int_{\partial \Omega_D} C_{11} v g_D \, ds +
\int_{\partial \Omega_N} v g_N \, ds \\
\bet
& = \int_{\partial \Omega_D} v \kappa \nabla_h u_h \cdot \bn \, ds +
\int_{\partial \Omega_D} \kappa v \, (r([u_h]) \nn \\
\bet
&\quad + l({\bm C}_{12}
\cdot [u_h]) + r_D(u_h)) \cdot \bn \, ds \\
\bet
& \quad -
\int_{\partial \Omega_D} \kappa v \, r_D(g_D)\cdot \bn \, ds +
\int_{\partial \Omega_D} C_{11} v (g_D - u_h) \, ds + \int_{\partial
\Omega_N} v g_N \, ds \\
\bet
& = \int_{\partial \Omega_D} v \kappa \nabla_h u_h \cdot \bn \, ds
- \int_{\Omega} \kappa r_D(v) \nn \\
\bet
&\quad \cdot (r([u_h]) + l({\bm C}_{12} \cdot
[u_h]) + r_D(u_h)) \, dx \\
\bet & \quad - \int_{\partial \Omega_D}\kappa v \, r_D(g_D)\cdot
\bn \, ds + \int_{\partial \Omega_D} C_{11} v (g_D - u_h) \, ds +
\int_{\partial \Omega_N} v g_N \, ds. \nonumber
\end{split}
\end{equation}
Therefore, we can rewrite (\ref{goveeq})~as
\begin{equation}
B_h^{LDG}(u_h, v) = L_h^{LDG}(v) \qquad \forall v \in V_h,
\end{equation}
where the bilinear form $B^{LDG}_h: V_h \times V_h \to \mathbb{R}$
is given by
%
\begin{align}
 B^{LDG}_h(u,v) & = \int_\Omega \nabla_h v\cdot (\kappa \nabla_h
u) \, dx - \int_{{\cal E}_i} ( [u] \cdot \{\kappa \nabla_h v\} +
\{\kappa \nabla_h u\} \cdot [v] ) \, ds\nonumber \\
\bet
 & \quad -
\int_{{\cal E}_i} ( {\bm C}_{12} \cdot [u] [\kappa \nabla_h v] +
[\kappa \nabla_h u]{\bm C}_{12} \cdot [v]) \, ds + \int_{{\cal
E}_i} C_{11}[u] \cdot
[v] \, ds \nonumber\\
\bet
 & \quad +
\int_{\Omega} \kappa (r([u]) + l({\bm C}_{12} \cdot [u]) + r_D(u)
) \cdot (r([v]) + l({\bm C}_{12} \cdot [v]) + r_D(v) ) \, dx \nonumber\\
\bet
& \quad - \int_{\partial \Omega_D} (\kappa \nabla_h u \cdot \bn v +
u \kappa \nabla_h v \cdot \bn) ds + \int_{\partial \Omega_D} C_{11}
u v \, ds \label{bilinear}
\end{align}
%
and the linear form $L^{LDG}_h: V_h \to \mathbb{R}$ is given by
\begin{align}
L^{LDG}_h(v) & = \int_\Omega f v \, dx
 - \int_{\partial \Omega_D} g_D ( \kappa \nabla_h v + r([v]) + l({\bm C}_{12}\cdot[v])) \cdot \bn ds\nonumber \\
\bet
 & \quad - \int_{\partial \Omega_D} \kappa v\, r_D(g_D) \cdot \bn \, ds + \int_{\partial \Omega_D} C_{11}
 g_D v \, ds
 + \int_{\partial \Omega_N} v\, g_N \, ds \quad
\forall v \in V_h . \label{LLDG}
\end{align}

It is straightforward to verify that the bilinear form
(\ref{bilinear}) is symmetric, i.e., $B_h(u,v) = B_h(v,u)$. Also,
the conservative form of the numerical fluxes, (\ref{qflux}) and
(\ref{uflux}), guarantees that the LDG scheme is conservative and
adjoint consistent \cite{arnold02unified}.

Unfortunately, when the scheme is implemented in multidimensions
on general triangular/tetrahedral meshes, the resulting
discretization is not compact in the sense that the equation
corresponding to a given degree of freedom may involve degrees of
freedom that belong to elements which are not immediate neighbors.
It turns out that these additional connections are due to the
volume term in (\ref{bilinear}) which involves products of the
lifting functions. Although the connectivity pattern between
elements depends on the choice of face switches in (\ref{switch}),
it is well known~\cite{sherwin06elliptic} that in multidimensions
this problem cannot be remedied by a more careful choice of the
face switches (\ref{switch}). This noncompactness of the LDG
scheme occurs also for quadrilateral/hexahedral discretizations.

\section{The CDG algorithm} \label{CDG}

The CDG algorithm is designed to be compact and, at the same time,
inherit all the attractive properties of the LDG algorithm. To
start with, we decompose the lifting operators introduced in
(\ref{liftings}) into facewise contributions. Thus, we consider
for all $e \in {\cal E}_i$, $r^e: [L^2(e)]^d \to \Sigma_h$, $l^e:
L^2(e) \to \Sigma_h$ and for each $e \in \partial \Omega_D$, $r_D: L^2(e)
\to \Sigma_h$, defined~as
\begin{equation}
\begin{array}{rcll}
\displaystyle \int_\Omega r^e(\phi) \cdot \btau \, dx & = &
\displaystyle - \int_{e} \phi \cdot \{ \btau \} \,
ds  \qquad & \forall \btau \in \Sigma_h, \\
\bet \displaystyle \int_\Omega l^e(q) \cdot \btau \, dx & = & -
\displaystyle \int_{e} q [ \btau] \, ds  & \forall \btau \in
\Sigma_h, \\
\bet
\displaystyle \int_\Omega r^e_D(q) \cdot \btau \, dx & = &
 \displaystyle - \int_{e} q \btau \cdot \bn \, ds
\qquad & \forall \btau \in \Sigma_h.
\end{array} \label{liftings_e}
\end{equation}
Clearly, we will have, for all $\phi \in [L^2({\cal E}_i)]^d$ and
all $q \in L^2({\cal E}_i)$,
\begin{equation}
r(\phi) = \sum_{e \in {\cal E}_i} r^e(\phi), \qquad l(q) = \sum_{e
\in {\cal E}_i} l^e(q), \qquad r_D(q) = \sum_{e \in \partial \Omega_D}
r^e_D(q).
\end{equation}

Now, we can define the CDG method. The numerical interelement
fluxes $(\hat{\bsigma}, \hat{u})$ for the CDG method are given~by
\begin{align}
\hat{\bsigma} & = \{ \bsigma^e_h \} - C_{11}[u_h] + {\bm C}_{12}
[\bsigma^e_h], \label{qfluxcdg} \\
\bet \hat{u} & = \{ u_h \} - {\bm C}_{12} \cdot [u_h]
\label{ufluxcdg}
\end{align}
for the interior faces, and
\begin{equation}
\begin{array}{r@{\,}c@{\,}l@{\,}r@{\,}cll}
 \hat{\bsigma} & = & \bsigma^e_h - C_{11}(u_h - g_D)\bn, \qquad
\hat{u} & = & g_D
\qquad & \mbox{on} \quad \partial \Omega_D, \\
\bet
 \hat{\bsigma} & = & g_N \bn, \hfill \hat{u} & = & u_h \qquad & \mbox{on}
\quad
\partial \Omega_N,
\end{array} \label{bfluxescdg}
\end{equation}
for the boundary faces. Here, $\bsigma^e_h$ is given~as
\begin{equation}
\bsigma^e_h = \kappa \nabla_h u_h + \bar{\bsigma}^e_h,
\end{equation}
where
\begin{equation}
\bar{\bsigma}^e_h = \kappa r^e([u_h]) + \kappa l^e({\bm
C}_{12}\cdot[u_h]) - \kappa r^e_D(g_D - u_h) \ .
\end{equation}
We note that the numerical flux, $\hat{u}$, is chosen as in the
LDG method. Therefore, (\ref{globalDGibpLDG}) and
(\ref{bsigmah})--(\ref{goveeq}) still apply for the CDG method,
and the only difference between the LDG and CDG methods is in the
evaluation of the terms involving $\hat{\sigma}$ in
(\ref{goveeq}), which in the CDG case is done according to
(\ref{qfluxcdg}) and (\ref{bfluxescdg}). Also, the coefficients
${\bm C}_{12}$ are given by expressions (\ref{switch})
and~(\ref{sumtoone}).

In order to compute the CDG numerical flux $\hat{\bsigma}$ on a
given face $e$, we need to evaluate first a stress field
$\bsigma^e_h$ associated with this face. This evaluation, however,
can be carried out efficiently due to the localized support of
$\bar{\bsigma}^e_h$. In particular, we note that when $e \in
\partial \Omega_N$, then $\bar{\bsigma}^e_h = {\bm 0}$. When $e \in
\partial \Omega_D$, we have $\bar{\bsigma}^e_h = \kappa r^e_D(g_D -
u_h)$, which has only a nonzero support on the element neighboring
face $e$. Finally, when $e \in {\cal E}_i$, then
$\bar{\bsigma}^e_h = \kappa r^e([u_h]) + \kappa l^e({\bm
C}_{12}\cdot[u_h])$. In this case, $\bar{\bsigma}^e_h$ is nonzero
only in one of the elements neighboring face $e$. The element in
which $\bar{\bsigma}^e_h$ is nonzero is determined by the choice
of switches for that face. In particular, using (\ref{switch}) and
(\ref{liftings_e}), it can be easily shown that if $S^{K^-}_{K^+}
= 1$ and $S^{K^+}_{K^-} = 0$, then $\bar{\bsigma}^e_h = {\bm 0}$
on $K^-$. Similarly, we will have $\bar{\bsigma}^e_h = {\bm 0}$ on
$K^+$ when $S^{K^-}_{K^+} = 0$ and $S^{K^+}_{K^-} = 1$.

\subsection{Primal form of the CDG algorithm}\label{sec3.1}

In order to obtain the primal form of the CDG method, we proceed
as before and start from (\ref{goveeq}). In this case, the terms
involving $\hat{\bsigma}$ become
\begin{equation}
\begin{split}
 & \int_{{\cal E}_i} \hat{\bsigma} \cdot [v] \, ds = \sum_{e \in {\cal E}_i}
 \int_{e} \hat{\bsigma} \cdot [v] \, ds \\
\bet
& =  \int_{{\cal E}_i}
( \{\kappa \nabla_h u_h \} + {\bm C}_{12} [ \kappa \nabla_h u_h ])
\cdot [v] \, ds + \sum_{e \in {\cal E}_i } \int_{e} (\{
\bar{\bsigma}^e \} + {\bm C}_{12} [ \bar{\bsigma}^e ]) \cdot [v] \,
ds \\
\bet
& \quad - \int_{{\cal E}_i} C_{11} [u]\cdot [v] \, ds \\
\bet
& = \int_{{\cal E}_i} ( \{\kappa \nabla_h u_h \} + {\bm C}_{12} [ \kappa
\nabla_h u_h ]) \cdot [v] \, ds \\
\bet
& \quad - \sum_{e \in {\cal E}_i } \int_{\Omega} \kappa (r^e([v]) + l^e({\bm C}_{12} \cdot [v]))\cdot
(r^e([u_h]) + l^e({\bm C}_{12} \cdot [u_h]) + r^e_D(u_h) ) \, dx \\
\bet
& \quad + \sum_{e \in {\cal E}_i }
\int_{\Omega} \kappa (r^e([v]) + l^e({\bm C}_{12} \cdot [v]))\cdot
r^e_D(g_D) \, dx - \int_{{\cal E}_i } C_{11} [u ] \cdot [v] \, ds
\nonumber
\end{split}
\end{equation}
and
\begin{equation}
\begin{split}
 & \int_{\partial
\Omega} v\, \hat{\bsigma} \cdot \bn \, ds = \sum_{e \in \partial
\Omega} \int_{e} v\, \hat{\bsigma} \cdot \bn \, ds \\
\bet
 & = \sum_{e \in \partial \Omega_D} \int_e v \bsigma^e_h \cdot \bn \,
ds - \int_{\partial \Omega_D} C_{11} v u_h \, ds + \int_{\partial
\Omega_D} C_{11} v g_D \, ds +
\int_{\partial \Omega_N} v g_N \, ds \\
\bet
 & = \int_{\partial \Omega_D} v \kappa \nabla_h u_h \cdot \bn \, ds +
\sum_{e \in \partial \Omega_D} \int_{e} \kappa v \, (r^e([u_h]) +
l^e({\bm C}_{12} \cdot [u_h]) + r^e_D(u_h)) \cdot \bn \, ds \\
\bet
 &\quad - \sum_{e \in \partial \Omega_D} \int_{e}
 \kappa v \, r^e_D(g_D)\cdot \bn \, ds
- \int_{\partial \Omega_D} C_{11} v u_h \, ds + \int_{\partial
\Omega_D} C_{11} v g_D \, ds + \int_{\partial \Omega_N} v g_N \, ds \\
\bet
& = \int_{\partial \Omega_D} v \kappa \nabla_h u_h \cdot \bn \, ds
- \sum_{e \in \partial \Omega_D} \int_{\Omega} \kappa r^e_D(v) \cdot
(r^e([u_h]) + l^e({\bm C}_{12} \cdot [u_h]) + r^e_D(u_h)) \, dx \\
\bet & \quad - \sum_{e \in \partial \Omega_D} \int_{e}\kappa v \,
r^e_D(g_D)\cdot \bn \, ds - \int_{\partial \Omega_D} C_{11} v u_h
\, ds + \int_{\partial \Omega_D} C_{11} v g_D \, ds +
\int_{\partial \Omega_N} v g_N \, ds. \nonumber
\end{split}
\end{equation}
Thus, for the CDG scheme, (\ref{goveeq}) can be written as
\begin{equation}
B_h^{CDG}(u_h, v) = L_h^{CDG}(v) \qquad \forall v \in V_h,
\end{equation}
where the bilinear form $B^{CDG}_h: V_h \times V_h \to \mathbb{R}$
is given~by
\begin{equation}
\begin{split}
 & B^{CDG}_h(u,v) = \int_\Omega
\nabla_h v\cdot (\kappa \nabla_h u) \, dx - \int_{{\cal E}_i} ( [u]
\cdot \{\kappa \nabla_h v\} +
\{\kappa \nabla_h u\} \cdot [v] ) \, ds \\
\bet
 & \quad -
\int_{{\cal E}_i} ( {\bm C}_{12} \cdot [u] [\kappa \nabla_h v] +
[\kappa \nabla_h u]{\bm C}_{12} \cdot [v]) \, ds \\
\bet
 & \quad +
\sum_{e \in ({\cal E}_i \bigcup \partial \Omega_D)} \int_{\Omega} \kappa (r^e([u]) + l^e({\bm
C}_{12} \cdot [u]) + r^e_D(u)
) \cdot (r^e([v]) + l^e({\bm C}_{12} \cdot [v]) + r^e_D(v) ) \, dx \\
\bet
& \quad - \int_{\partial \Omega_D} (\kappa \nabla_h u \cdot \bn v +
u \kappa \nabla_h v \cdot \bn) ds + \int_{{\cal E}_i} C_{11}[u]
\cdot [v] \, ds + \int_{\partial \Omega_D} C_{11} u v \, ds
\label{bilinearcdg}
\end{split}
\end{equation}
and the linear form $L^{CDG} : V_h \to \mathbb{R}$ is given by 
\begin{align}
L^{CDG}_h(v) & = \int_\Omega f v \, dx
 - \int_{\partial \Omega_D} g_D  \kappa \nabla_h v  \cdot \bn ds\nonumber \\
\bet
 & \quad - \int_{\partial \Omega_D} \kappa v\, r_D(g_D) \cdot \bn \, ds + \int_{\partial \Omega_D} C_{11}
 g_D v \, ds
 + \int_{\partial \Omega_N} v\, g_N \, ds \quad
\forall v \in V_h . \label{LCDG}
\end{align}

The CDG method is symmetric, i.e.,\ $B^{CDG}_h(u,v) =
B^{CDG}_h(v,u)$, and retains all the attractive properties of the
LDG algorithm such as consistency and adjoint consistency.

\subsection{Error estimates}\label{sec3.2}

We observe that the only difference between the LDG and CDG
schemes is the stabilizing term involving the products of the
lifting functions. In the LDG scheme, we have
\begin{equation}
\begin{split}
& \int_{\Omega} \kappa (r([u]) + l({\bm C}_{12} \cdot [u]) + r_D(u) )
\cdot (r([v]) + l({\bm C}_{12} \cdot [v]) + r_D(v) ) \, dx, \\
\bet & \sum_{e \in {\cal E}_i} \sum_{f \in {\cal E}_i}
\int_{\Omega} \kappa (r^e([u]) + l^e({\bm C}_{12} \cdot [u]) +
r^e_D(u) ) \cdot (r^f([v]) + l^f({\bm C}_{12} \cdot [v]) +
r^f_D(v) ) \, dx,
\label{LDGst}
\end{split}
\end{equation}
whereas in the CDG scheme, we have
\begin{equation}
\begin{split}
& \sum_{e \in {\cal E}_i} \int_{\Omega} \kappa (r^e([u]) + l^e({\bm
C}_{12} \cdot [u]) + r^e_D(u) ) \cdot (r^e([v]) + l^e({\bm C}_{12}
\cdot [v]) + r^e_D(v) ) \, dx, \\
\bet
& \sum_{e \in {\cal E}_i} \sum_{f \in {\cal E}_i} \delta_{ef} \int_{\Omega} \kappa (r^e([u]) + l^e({\bm
C}_{12} \cdot [u]) + r^e_D(u) ) \cdot (r^f([v]) + l^f({\bm C}_{12}
\cdot [v]) + r^f_D(v) ) \, dx , \label{CDGst}
\end{split}
\end{equation}
where $\delta_{ef}$ is the Kronecker delta. Thus, we see that the
CDG scheme can be regarded as the LDG algorithm with some terms
turned off. We also note that the turned-off terms in the LDG
algorithm are indefinite and hence are not guaranteed to
contribute to the method's stability. The effect of using lifting
functions in the CDG method which are associated with individual
faces is to eliminate connectivities between nonneighboring
elements. We note that an analogous approach was adopted
in~\cite{bassi97br2,brezzi00br2} to render the BR2 scheme compact.

It turns out that the proofs of coercivity and boundedness for the
LDG method presented in~\cite{arnold02unified} can be used here
without change. This leads to optimal a priori estimates for the
CDG method,
\begin{equation}
||| u-u_h ||| \le C h^p |u|_{p+1,\Omega}
\end{equation}
and
\begin{equation}
|| u-u_h ||_{0,\Omega} \le C h^{p+1} |u|_{p+1,\Omega} \ .
\end{equation}
Here, the norm $||| \cdot |||$ is given~by
\begin{equation}
||| v |||^2 = \sum_{K \in {\cal T}_h}| v|^2_{1,K} + \sum_{e \in
{\cal E}_i} || r_e([v])||^2_{0,\Omega} + \sum_{e \in \partial
{\Omega_D}} || r_D(v)||^2_{0,\Omega} \ .
\end{equation}

The above estimates require that the stabilization parameter
$C_{11}$ in (\ref{qfluxcdg}) is taken to be of order ${\cal
O}(h^{-1})$, where $h$ is the characteristic mesh size (see
also~\cite{castillo00ldg}). We note that for $C_{11}$ of order
${\cal O}(1)$, only suboptimal convergence is demonstrated, but in
practical computations, optimal results are also observed. We also
point out that for general discretizations, the piecewise constant
approximation $p=0$ does not lead to a consistent discretization.
This is in common with other DG schemes such as the LDG or
the~BR2.

\section{Stabilization}\label{STAB}

The above a priori error estimates are applicable to both the CDG
and LDG algorithms. It turns out that, for the LDG algorithm, one
can set $C_{11}=0$ for all the internal interfaces, provided the
switches in (\ref{switch}) are chosen following a simple rule.
That is, if the switches for each simplex element $K$ satisfy that
 \begin{equation}
 \sum_{e\in \partial K}
S_K^{K'} < d + 1,
\label{ruleswitch}
\end{equation}
where $d$ is the problem dimension, then the scheme shows no
degradation in performance and becomes extremely simple. This
result was proven in~\cite{MDL}. In this case, the numerical flux
$\hat{u}$ on a given internal face is taken to be the value of
$u_h$ on one of the neighboring elements, while the numerical flux
$\hat{\bsigma}$ is taken to be the value of $\bsigma_h$ on the
other neighboring element. The element used to calculate either
$\hat{u}$ or $\hat{\bsigma}$ is determined by the value of
switches on that face. The rule (\ref{ruleswitch}) guarantees
that, when calculating the numerical fluxes on each face, the
value of the solution on each element will be used, at least once,
to set $\hat{u}$ on the element boundary, and, at least once, to
set $\hat{\bsigma}$ on the element boundary.

Clearly, there is plenty of flexibility in choosing appropriate
values for switches which satisfy the rule (\ref{ruleswitch}); see
\cite{cockburncartesian}, for instance. Thus, provided that the
rule (\ref{ruleswitch}) is satisfied, the LDG scheme converges at
the optimal rate without the need for explicit stabilization.

\subsection{Null-space dimension}\label{sec4.1}

We have found that while the rule (\ref{ruleswitch}) is essential
in ensuring that the solution is unique for the LDG method, this
requirement is not necessary for the CDG method. That is, for the
CDG method we are able to set $C_{11} =0$ for all the internal
faces and use any combination of switches with the only constraint
given by~(\ref{sumtoone}).

In order to illustrate this point, we adopt the two-dimensional
test problem presented in \cite{sherwin06elliptic}. We consider a
square domain with periodic boundary conditions imposed on all
sides. We perform a regular subdivision into four squares and then
subdivide each square into two triangles. We look at
approximations ranging from $p=1$ to $p=7$ and nodal basis
functions with equally spaced nodes. We discretize the Laplacian
operator using the CDG and the LDG algorithms with the parameter
$C_{11}$ set to zero and calculate the dimension of the null-space
of the resulting matrix.

We consider two different switches for both the LDG and CDG
algorithms. The so-called {\em consistent switch} satisfies
(\ref{ruleswitch}), and here it is chosen using a procedure
analogous to that presented in~\cite{cockburncartesian,
sherwin06elliptic}. We also consider the {\em natural switch},
which is based on element numbering and sets $S_{K^+}^{K^-} = 0$
if the element number $K^+$ is less than the element number $K^-$,
and to $1$ otherwise. This switch was first introduced
in~\cite{consistent} in the context of interior point methods for
elliptic problems.

Because of the periodic boundary conditions, any solution will be
undetermined up to a constant, and as a consequence, we expect a
singular matrix with a null-space of dimension one. The computed
dimension of the null-space for the different schemes, polynomial
order interpolations, and switches is presented in
Table~\ref{tab:nullspace}. We note that while the LDG scheme gives
the desired null-space dimension of one when the consistent switch
is employed, the null-space dimension grows with increasing $p$,
when the natural switch is employed. This same result was reported
in \cite{sherwin06elliptic}. On the other hand, the CDG scheme
always gives the desired one-dimensional null-space for all $p$
and for any switch choice.

We note that the natural switch has some computational advantages
when computing the ILU(0) factorization of the system
matrix~\cite{implicit}. If $S_{K^+}^{K^-} = 0$ when $K^+ < K^-$,
the lower triangular blocks in the matrix have only a few nonzero
rows, and no additional fill-in is introduced during the
factorization phase. On the other hand, for an arbitrary switch
choice, some lower triangular blocks will have nonzero columns
that will render the blocks completely full after factorization.
This effect is described in more detail in \cite{implicit}, where
the CDG method is used to discretize convective-diffusive systems
which are solved using a preconditioned Krylov solver.

\begin{table}
\begin{center}
Nullspace dimension \\
\begin{tabular}{ll|ccccccc}
\hline
Polynomial order $p$ & & $1$ & $2$ & $3$ & $4$ & $5$ & $6$ & $7$ \\
\hline
Consistent switch & CDG & $1$ & $1$ & $1$ & $1$ & $1$ & $1$ & $1$ \\
 & LDG & $1$ & $1$ & $1$ & $1$ & $1$ & $1$ & $1$ \\
Natural switch & CDG & $1$ & $1$ & $1$ & $1$ & $1$ & $1$ & $1$ \\
 & LDG & $3$ & $4$ & $5$ & $6$ & $7$ & $8$ & $9$ \\
\hline
\end{tabular}
 \caption{Nullspace dimensions for the CDG/LDG schemes using the two
different switches. The problem is expected to have a
one-dimensional null-space, but with the (inconsistent) natural
switch the LDG scheme gives spurious modes and a null-space that
grows with~$p$.} \label{tab:nullspace}
\end{center}
\end{table}

\section{Implementation}\label{IMPLEMENT}

Since the main motivation for developing the CDG algorithm is to
obtain a computationally more efficient method, we next discuss
some practical implementation issues.

\begin{figure}
\begin{minipage}{.45\textwidth}
\begin{minipage}{\textwidth}
\hfill
\includegraphics[width=.45\textwidth]{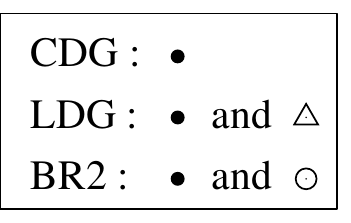}
\end{minipage}
\begin{minipage}{\textwidth}
\includegraphics[width=.85\textwidth]{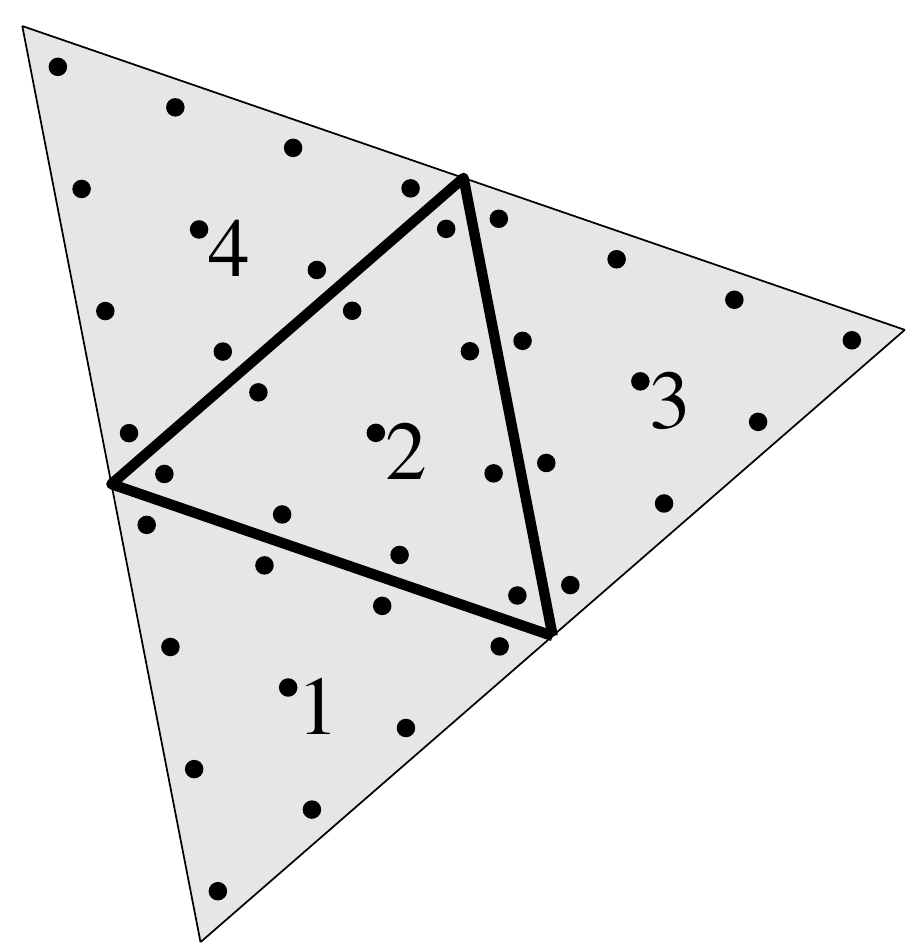}
\end{minipage}
\end{minipage}
\begin{minipage}{.55\textwidth}
\includegraphics[width=\textwidth]{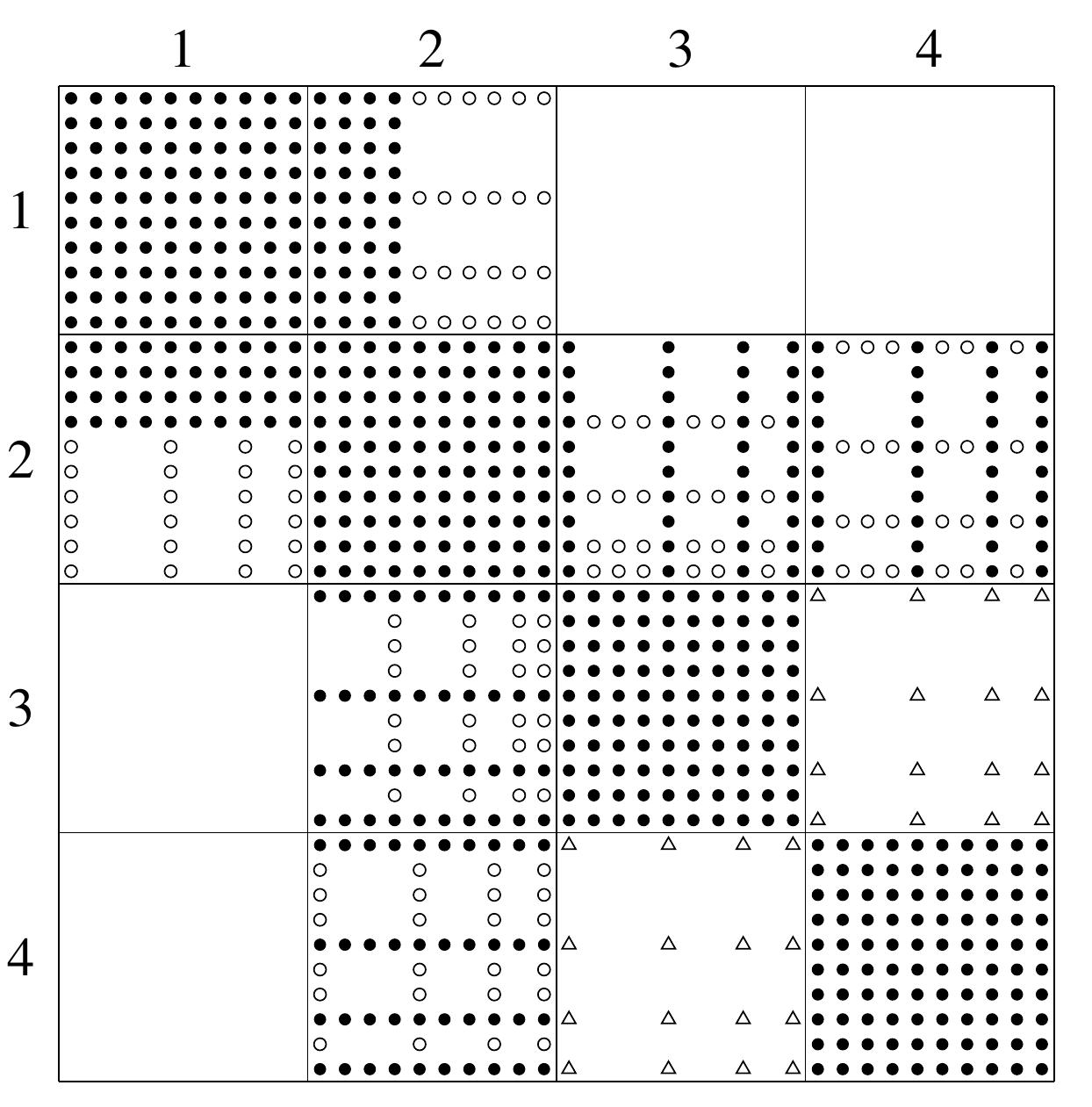}
\end{minipage}
\caption{The sparsity structure for four triangles
with $p=3$ (left plot). The CDG and the BR$2$ scheme are both
compact in the sense that they connect only neighboring triangles;
however, BR$2$ introduces more nonzeros. The LDG scheme is
noncompact and gives connections between some nonneighboring
triangles ($3$ and $4$).} \label{fig:sparsity}
\end{figure}

\subsection{Sparsity patterns}\label{sec5.1}

We start by discussing the sparsity pattern of the CDG
method and compare it with that of the LDG and BR2 methods. We
assume throughout that nodal bases~\cite{nodal} are used to span
the approximating and weighting Galerkin spaces. For illustration
purposes, we consider the triangular mesh shown in
Figure~\ref{fig:sparsity}, consisting of four elements and a
finite element space of piecewise polynomials of degree $p=3$ on
each element. The total number of degrees of freedom is $60$,
corresponding to $15$ degrees of freedom per element. The sparsity
patterns corresponding to the CDG, LDG, and BR2 methods are also
shown in Figure~\ref{fig:sparsity}. We note that the sparsity
pattern of the IP method is identical to that of the BR2 method,
and therefore the same remarks~apply.

As is well known, the LDG scheme introduces connections between
degrees of freedom in nonneighboring elements. In this example,
some degrees of freedom in element $3$ are connected to degrees of
freedom in element $4$. These connections are caused by the
stabilization term (\ref{LDGst}), which involves the product of
global lifting functions. We note that these nonlocal
connectivities also occur for quadrilateral discretizations and
cannot be avoided by a more careful renumbering of the elements
and/or internal interfaces~\cite{sherwin06elliptic}.

Of the three schemes, the CDG method produces the smallest number
of nonzero entries in the matrix. In fact, any nonzero entry in
the CDG matrix is also a nonzero entry in the matrices generated
by the other two schemes. The BR2 scheme is compact but connects
the face nodes of each element with all the nodes of the
neighboring element sharing that face. On the other hand, the CDG
scheme connects only the nodes of those faces for which the switch
is one, to the interior nodes of the neighboring element sharing
that face.

\subsection{Storage requirements}\label{sec5.2}

In order to quantify the matrix storage requirements for
the three schemes, we consider a simplex element in $d$ dimensions
having $d+1$ distinct neighboring elements. For polynomial basis
function of degree $p$, the number of degrees of freedom per
element is given by $S=\binom{p+d}{d}$, and the number of degrees
of freedom along each element face is given by
$S_e=\binom{p+d-1}{d-1}$. Using this notation, we can obtain
expressions for the number of nonzero matrix entries per interior
element.

For the CDG scheme we have one diagonal block with $S^2$ entries
and $d+1$ off-diagonal blocks with $S_eS$ entries. Since the
scheme connects some element face nodes to all the nodes of the
neighboring element sharing that face, we have
\begin{align*}
 M_\mathrm{CDG} = S^2+(d+1)S_eS.
\end{align*}

For the LDG scheme, the pattern is the same as for the CDG
algorithm plus the additional nonlocal connectivities. Each such
connectivity involves $S_e^2$ entries since the scheme connects
face nodes to nonneighboring face nodes. The number of nonlocal
connections $\alpha$ depends on the mesh and the switch, but on
average, we have that in one dimension the switch can be chosen
such that $\alpha=0$, and our experiments indicate that
$\alpha\approx 1$ for $d=2$ and $\alpha\approx 2$ for $d=3$. The
total number of nonzeros is then
\begin{align*}
 M_\mathrm{LDG} = S^2+(d+1)S_eS+\alpha S_e^2.
\end{align*}

Finally, for the BR2 (and also the IP) scheme, the pattern is the
same as with the CDG scheme, but with the additional connections
caused by the fact that all the face nodes connect to all the
interior nodes in the neighboring elements. This results in
$S_eS+S_e(S-Se)$ entries per block, giving a total number of
nonzeros of
\begin{align*}
 M_\mathrm{BR2} = S^2+(d+1)(2S-S_e)S_e.
\end{align*}

The memory requirements for $d=1,2,3$ and $p=1,\ldots,5$ are shown
in Table~\ref{tab:memory}. We note that the CDG method has the
lowest memory requirements. For instance, in three dimensions with
polynomials of degree $p=4$, the additional storage requirements
of the LDG and BR2 methods are 14\% and 36\%, respectively.

\begin{table}
\begin{center}
\begin{tabular}{cl|rrrrr}
\hline
 Dim & Scheme & \multicolumn{1}{c}{$p=1$} & \multicolumn{1}{c}{$p=2$} &
 \multicolumn{1}{c}{$p=3$} & \multicolumn{1}{c}{$p=4$} &
 \multicolumn{1}{c}{$p=5$} \\ \hline
 1  & CDG & 8 & 15 & 24 & 35 & 48 \\
 & LDG & 8 & 15 & 24 & 35 & 48 \\
 & BR2 & 10 & 19 & 30 & 43 & 58 \\
 & & & & & \\
 2  & CDG & 27 & 90 & 220 & 450 & 819 \\
 & LDG & 31 & 99 & 236 & 475 & 855 \\
 & BR2 & 33 & 117 & 292 & 600 & 1089 \\
 & & & & & \\
 3  & CDG & 64 & 340 & 1200 & 3325 & 7840 \\
 & LDG & 82 & 412 & 1400 & 3775 & 8722 \\
 & BR2 & 76 & 436 & 1600 & 4525 & 10780 \\
\hline
\end{tabular}
 \caption{Memory requirements per interior simplex
element for the CDG, LDG, and BR$2$ schemes. The case $p=3$ in two
dimensions is illustrated in Figure~$\ref{fig:sparsity}.$ The LDG
scheme is assumed to have $\alpha=0,1,2$ noncompact neighbors in
one, two, and three dimensions, respectively.} \label{tab:memory}
\end{center}
\end{table}

Finally, we note that the CDG sparsity pattern is such that in
addition to having fewer nonzero entries, the entire matrix can be
stored using simple blockwise dense arrays. In particular, for a
problem involving $T$ elements, we can use an $S\times S \times T$
dense array for the diagonal blocks, and an $S \times S_e \times
(d+1) \times T$ dense array for the off-diagonal blocks. This
representation is not only simple and compact, it also makes it
straightforward to apply high-performance libraries such as the
BLAS routines \cite{BLAS} for basic matrix operations.

A similar storage format is harder to define for the LDG scheme,
because of the noncompactness and the somewhat complex pattern in
which these additional blocks appear. For the BR2 scheme, while it
is compact, and in principle one could use a storage scheme
similar to that of the CDG method, the sparsity pattern of the
off-diagonal blocks is nonrectangular, and therefore any dense
storage strategy would require, at least, an additional array.

\section{Numerical results}\label{NUMERICAL}

In this section, we present some numerical experiments to assess
the accuracy and behavior of the CDG algorithm. We consider a
two-dimensional model problem. The domain $\Omega$ is the unit
square $[0,1]\times [0,1]$. Dirichlet conditions are imposed at
all the boundaries ($\partial \Omega_D =\partial \Omega$), and we
choose the analytical solution
\begin{align}
 u(x,y) = \mathrm{exp} \left[ \alpha \sin (ax+by) + \beta \cos (cx+dy)
 \right] \label{solution}
\end{align}
with numerical parameters $\alpha=0.1, \beta=0.3, a=5.1, b=-6.2,
c=4.3, d=3.4$. We then solve the model Poisson problem (\ref{probdef}) with
the parameter $\kappa=1$ and with the Dirichlet boundary conditions
$g_D(x,y)=u(x,y)|_{\partial \Omega_D}$. The source term, $f(x,y)$,
is obtained by analytical differentiation of~(\ref{solution}).

We consider triangular meshes obtained by splitting a regular $n
\times n$ Cartesian grid into a total of $2n^2$ triangles, giving
uniform element sizes of $h=1/n$. On these meshes, we consider
solutions of polynomial degree $p$ represented using a nodal basis
within each triangle, with the nodes uniformly distributed. We use
five different meshes, $n=2,4,8,16,32$, and five polynomial
degrees, $p=1$ to $p=5$.

\subsection{\boldmath Effect of the stabilization parameter $C_{11}$}
\label{sec6.1}

In order to assess the effect of the stabilization parameter, we
discretize the Poisson equation (\ref{probdef}) in two dimensions
and solve for the numerical solution $u_h$ using different values
of the stabilization parameter $C_{11}$. The resulting equation
system is solved using a preconditioned iterative solver
\cite{implicit}. We then compute the $L_2$ error $||u -
u_h||_{0,\Omega}$. The computed $L_2$ error, $||u -
u_h||_{0,\Omega}$, is shown in Table~\ref{tab:L2consistent} for
the different values of $p$ and $n$, and for $C_{11}=0,1,$ and
$10$, using the consistent switch. The same results are reported
for the natural switch in Table~\ref{tab:L2natural}.

We note that the accuracy is only weakly dependent on the value of
$C_{11}$. The only noticeable differences are for the
underresolved cases ($p=1,2$ and $n=2$) when using a large amount
of stabilization, $C_{11}=10$. We obtain the optimal convergence
rate of $p+1$ for all cases. Using the natural switch, instead of
the consistent one, makes the errors somewhat larger, but on
average only by 11\% and, in the worst case, only by~42\%.

\begin{table}
\begin{center}
\begin{tabular}{lr|rrrrr|c}
 \hline
$p$ & $C_{11}$ & \multicolumn{1}{c}{$n=2$} &
\multicolumn{1}{c}{$n=4$} &
 \multicolumn{1}{c}{$n=8$} & \multicolumn{1}{c}{$n=16$} &
 \multicolumn{1}{c|}{$n=32$} & Rate \\
  \hline
\multicolumn{1}{c}{} & \multicolumn{1}{c}{} &
\multicolumn{1}{|c}{} & \multicolumn{1}{c}{} &
\multicolumn{1}{c}{} &
\multicolumn{1}{c}{} & \multicolumn{1}{c|}{}& \multicolumn{1}{c}{}
\\[-7pt]
1 & 0 & $4.55\cdot 10^{-2}$ & $1.52\cdot 10^{-2}$ & $4.63\cdot 10^{-3}$ & $1.26\cdot 10^{-3}$ & $3.27\cdot 10^{-4}$ \ & 1.9 \\
 & 1 & $4.55\cdot 10^{-2}$ & $1.49\cdot 10^{-2}$ & $4.56\cdot 10^{-3}$ & $1.25\cdot 10^{-3}$ & $3.26\cdot 10^{-4}$ \ & 1.9 \\
 & 10 & $2.20\cdot 10^{-0}$ & $2.07\cdot 10^{-2}$ & $4.24\cdot 10^{-3}$ & $1.16\cdot 10^{-3}$ & $3.13\cdot 10^{-4}$ \ & 1.9 \\
 & & & & & & \\
2 & 0 & $9.00\cdot 10^{-3}$ & $1.80\cdot 10^{-3}$ & $2.56\cdot 10^{-4}$ & $3.36\cdot 10^{-5}$ & $4.29\cdot 10^{-6}$ \ & 3.0 \\
 & 1 & $9.10\cdot 10^{-3}$ & $1.80\cdot 10^{-3}$ & $2.56\cdot 10^{-4}$ & $3.36\cdot 10^{-5}$ & $4.29\cdot 10^{-6}$ \ & 3.0 \\
 & 10 & $2.89\cdot 10^{-2}$ & $2.01\cdot 10^{-3}$ & $2.62\cdot 10^{-4}$ & $3.38\cdot 10^{-5}$ & $4.30\cdot 10^{-6}$ \ & 3.0 \\
 & & & & & & \\
3 & 0 & $2.61\cdot 10^{-3}$ & $2.44\cdot 10^{-4}$ & $1.72\cdot 10^{-5}$ & $1.11\cdot 10^{-6}$ & $7.04\cdot 10^{-8}$ \ & 4.0 \\
 & 1 & $2.63\cdot 10^{-3}$ & $2.44\cdot 10^{-4}$ & $1.72\cdot 10^{-5}$ & $1.11\cdot 10^{-6}$ & $7.04\cdot 10^{-8}$ \ & 4.0 \\
 & 10 & $4.16\cdot 10^{-3}$ & $2.59\cdot 10^{-4}$ & $1.73\cdot 10^{-5}$ & $1.11\cdot 10^{-6}$ & $7.03\cdot 10^{-8}$ \ & 4.0 \\
 & & & & & & \\
4 & 0 & $1.09\cdot 10^{-3}$ & $4.52\cdot 10^{-5}$ & $1.57\cdot 10^{-6}$ & $5.14\cdot 10^{-8}$ & $1.64\cdot 10^{-9}$ \ & 5.0 \\
 & 1 & $1.09\cdot 10^{-3}$ & $4.54\cdot 10^{-5}$ & $1.57\cdot 10^{-6}$ & $5.15\cdot 10^{-8}$ & $1.64\cdot 10^{-9}$ \ & 5.0 \\
 & 10 & $1.19\cdot 10^{-3}$ & $4.77\cdot 10^{-5}$ & $1.60\cdot 10^{-6}$ & $5.16\cdot 10^{-8}$ & $1.64\cdot 10^{-9}$ \ & 5.0 \\
 & & & & & & \\
5 & 0 & $3.73\cdot 10^{-4}$ & $9.31\cdot 10^{-6}$ & $1.76\cdot 10^{-7}$ & $2.83\cdot 10^{-9}$ & $4.47\cdot 10^{-11}$ & 6.0 \\
 & 1 & $3.75\cdot 10^{-4}$ & $9.32\cdot 10^{-6}$ & $1.76\cdot 10^{-7}$ & $2.83\cdot 10^{-9}$ & $4.47\cdot 10^{-11}$ & 6.0 \\
 & 10 & $4.07\cdot 10^{-4}$ & $9.52\cdot 10^{-6}$ & $1.77\cdot 10^{-7}$ & $2.84\cdot 10^{-9}$ & $4.47\cdot 10^{-11}$ & 6.0 \\
\hline
\end{tabular}
 \caption{$L_2$ errors in the solution for the model Poisson problem, for various polynomial degrees $p$, mesh sizes $n$, and $C_{11}$ values. The consistent switch is used. The convergence rate is calculated based on the two finest meshes.}
\label{tab:L2consistent}
\end{center}
\end{table}

\begin{table}
\begin{center}
\begin{tabular}{lr|rrrrr|c}
\hline
$p$ & $C_{11}$ & \multicolumn{1}{c}{$n=2$} & \multicolumn{1}{c}{$n=4$} &
 \multicolumn{1}{c}{$n=8$} & \multicolumn{1}{c}{$n=16$} &
 \multicolumn{1}{c}{$n=32$} & Rate \\
 \hline
\multicolumn{1}{c}{} & \multicolumn{1}{c}{} &
\multicolumn{1}{|c}{} & \multicolumn{1}{c}{} &
\multicolumn{1}{c}{} & \multicolumn{1}{c}{} &
\multicolumn{1}{c|}{}& \multicolumn{1}{c}{}
\\[-7pt]
1 & 0 & $3.72\cdot 10^{-2}$ & $1.61\cdot 10^{-2}$ & $4.71\cdot 10^{-3}$ & $1.30\cdot 10^{-3}$ & $3.39\cdot 10^{-4}$ \ & 1.9 \\
 & 1 & $3.83\cdot 10^{-2}$ & $1.50\cdot 10^{-2}$ & $4.70\cdot 10^{-3}$ & $1.32\cdot 10^{-3}$ & $3.38\cdot 10^{-4}$ \ & 2.0 \\
 & 10 & $2.33\cdot 10^{-1}$ & $3.40\cdot 10^{-2}$ & $4.64\cdot 10^{-3}$ & $1.25\cdot 10^{-3}$ & $3.31\cdot 10^{-4}$ \ & 1.9 \\
 & & & & & & \\
2 & 0 & $1.28\cdot 10^{-2}$ & $1.96\cdot 10^{-3}$ & $3.03\cdot 10^{-4}$ & $3.98\cdot 10^{-5}$ & $5.04\cdot 10^{-6}$ \ & 3.0 \\
 & 1 & $1.18\cdot 10^{-2}$ & $2.07\cdot 10^{-3}$ & $2.88\cdot 10^{-4}$ & $4.01\cdot 10^{-5}$ & $5.02\cdot 10^{-6}$ \ & 3.0 \\
 & 10 & $3.24\cdot 10^{-2}$ & $3.00\cdot 10^{-3}$ & $3.37\cdot 10^{-4}$ & $4.05\cdot 10^{-5}$ & $5.16\cdot 10^{-6}$ \ & 3.0 \\
 & & & & & & \\
3 & 0 & $3.03\cdot 10^{-3}$ & $2.68\cdot 10^{-4}$ & $2.01\cdot 10^{-5}$ & $1.33\cdot 10^{-6}$ & $8.63\cdot 10^{-8}$ \ & 4.0 \\
 & 1 & $3.25\cdot 10^{-3}$ & $2.74\cdot 10^{-4}$ & $2.05\cdot 10^{-5}$ & $1.33\cdot 10^{-6}$ & $8.61\cdot 10^{-8}$ \ & 3.9 \\
 & 10 & $1.84\cdot 10^{-2}$ & $3.56\cdot 10^{-4}$ & $2.28\cdot 10^{-5}$ & $1.38\cdot 10^{-6}$ & $8.79\cdot 10^{-8}$ \ & 4.0 \\
 & & & & & & \\
4 & 0 & $9.67\cdot 10^{-4}$ & $5.15\cdot 10^{-5}$ & $1.82\cdot 10^{-6}$ & $5.86\cdot 10^{-8}$ & $1.87\cdot 10^{-9}$ \ & 5.0 \\
 & 1 & $1.33\cdot 10^{-3}$ & $5.24\cdot 10^{-5}$ & $1.81\cdot 10^{-6}$ & $5.90\cdot 10^{-8}$ & $1.88\cdot 10^{-9}$ \ & 5.0 \\
 & 10 & $1.98\cdot 10^{-3}$ & $6.30\cdot 10^{-5}$ & $1.95\cdot 10^{-6}$ & $6.12\cdot 10^{-8}$ & $1.90\cdot 10^{-9}$ \ & 5.0 \\
 & & & & & & \\
5 & 0 & $3.98\cdot 10^{-4}$ & $1.01\cdot 10^{-5}$ & $1.85\cdot 10^{-7}$ & $3.07\cdot 10^{-9}$ & $4.83\cdot 10^{-11}$ & 6.0 \\
 & 1 & $3.84\cdot 10^{-4}$ & $1.02\cdot 10^{-5}$ & $1.88\cdot 10^{-7}$ & $3.06\cdot 10^{-9}$ & $4.85\cdot 10^{-11}$ & 6.0 \\
 & 10 & $4.97\cdot 10^{-4}$ & $1.15\cdot 10^{-5}$ & $1.98\cdot 10^{-7}$ & $3.11\cdot 10^{-9}$ & $4.88\cdot 10^{-11}$ & 6.0 \\
\hline
\end{tabular}
 \caption{$L_2$ errors in the solution for the model Poisson problem, for various polynomial degrees
 $p$, mesh sizes $n$, and $C_{11}$ values. The natural switch is used. }
\label{tab:L2natural}
\end{center}
\end{table}

Table~\ref{tab:L2grad} shows the errors and the convergence rates
for the gradient of the solution using the CDG method with $C_{11}
= 0$. In particular, we calculate the seminorm $(\sum_{K \in {\cal
T}_h} |u-u_h|^2_{1,K})^{1/2}$. We observe optimal convergence at
the expected rate of~$p$.

\begin{table}
\begin{center}
\begin{tabular}{l|rrrrr|c}
\hline
$p$ & \multicolumn{1}{c}{$n=2$} & \multicolumn{1}{c}{$n=4$} &
 \multicolumn{1}{c}{$n=8$} & \multicolumn{1}{c}{$n=16$} &
 \multicolumn{1}{c|}{$n=32$} & Rate \\ \hline
 \multicolumn{1}{c|}{} & \multicolumn{1}{c}{} & \multicolumn{1}{c}{} &
 \multicolumn{1}{c}{} & \multicolumn{1}{c}{} &
 \multicolumn{1}{c}{} & \multicolumn{1}{|c}{} \\[-7pt]
 1 & $1.80\cdot 10^{-0}$ & $6.09\cdot 10^{-1}$ & $3.05\cdot 10^{-1}$ & $1.54\cdot 10^{-1}$ & $7.75\cdot 10^{-2}$ \ & 1.0 \\
 2 & $7.40\cdot 10^{-1}$ & $1.57\cdot 10^{-1}$ & $3.73\cdot 10^{-2}$ & $9.20\cdot 10^{-3}$ & $2.28\cdot 10^{-3}$ \ & 2.0 \\
 3 & $2.57\cdot 10^{-1}$ & $3.01\cdot 10^{-2}$ & $3.63\cdot 10^{-3}$ & $4.37\cdot 10^{-4}$ & $5.36\cdot 10^{-5}$ \ & 3.0 \\
 4 & $9.53\cdot 10^{-2}$ & $5.96\cdot 10^{-3}$ & $3.61\cdot 10^{-4}$ & $2.18\cdot 10^{-5}$ & $1.32\cdot 10^{-6}$ \ & 4.0 \\
 5 & $5.42\cdot 10^{-2}$ & $1.33\cdot 10^{-3}$ & $3.67\cdot 10^{-5}$ & $1.04\cdot 10^{-6}$ & $3.11\cdot 10^{-8}$ \ & 5.0 \\
\hline
\end{tabular}
\caption{The errors in the gradient for the CDG scheme with
consistent switch and $C_{11}=0$.}
\label{tab:L2grad}
\end{center}
\end{table}

\subsection{Comparison with the LDG and BR2 schemes}\label{sec6.2}

Here, we discretize the equations using the CDG, LDG, and BR2
schemes. For the CDG and the LDG methods, we use the consistent
switch and set $C_{11}=0$, except at the Dirichlet boundaries,
where $C_{11}=1$. The lifting parameter in the BR2 scheme is
$\eta=3$, which is the value required for
stability~\cite{brezzi00br2}.

\begin{figure}
\begin{center}
\includegraphics[width=.7\textwidth]{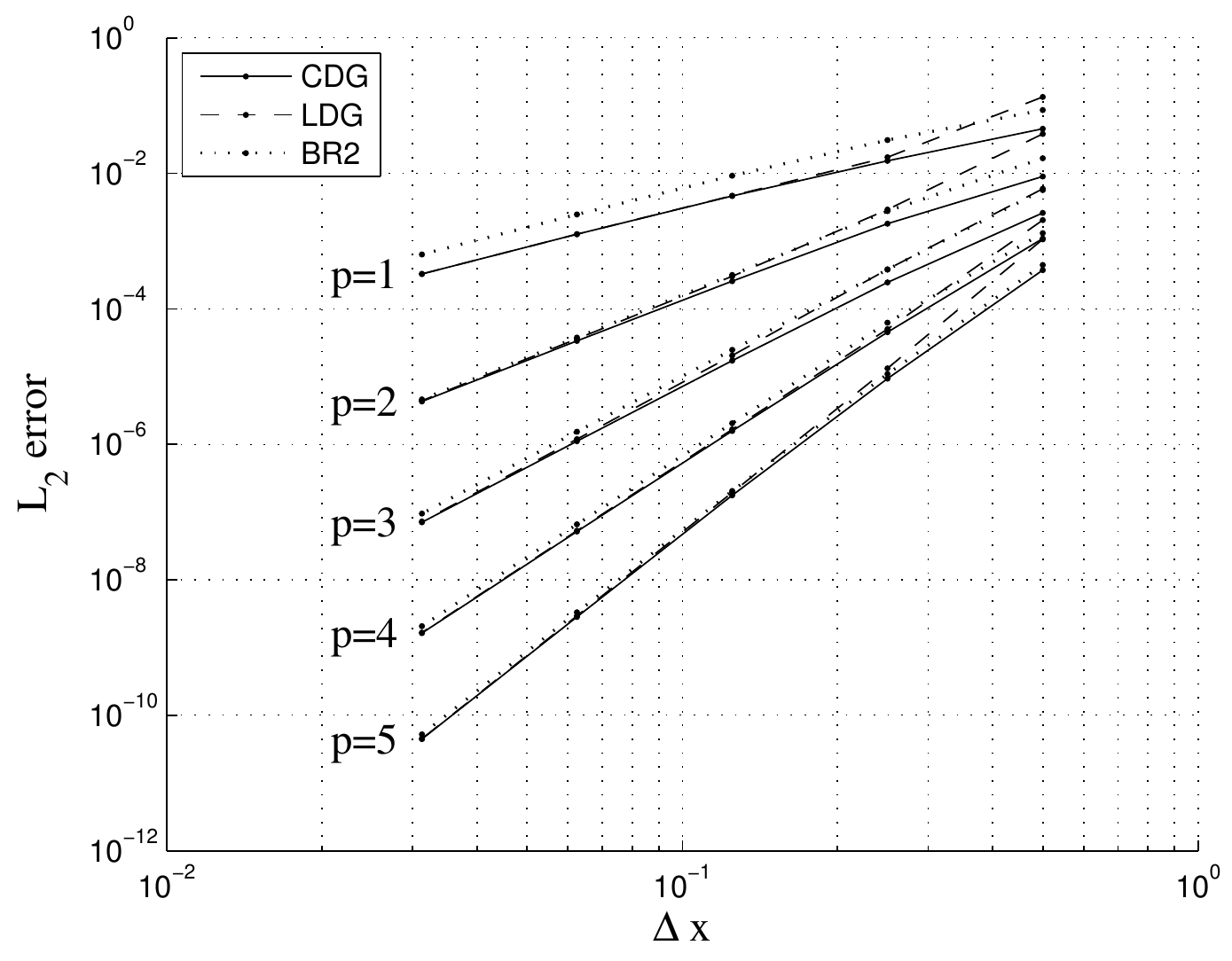}
\end{center}
\caption{$L_2$ errors in the solution for the model Poisson
problem; see Table~$\ref{tab:L2sol}$ for detailed values and
convergence rates.} \label{fig:L2sol}
\end{figure}

The accuracy results for the CDG, LDG, and BR2 schemes are shown
in Figure~\ref{fig:L2sol}, with details in Table~\ref{tab:L2sol}.
We note that the CDG scheme is the most accurate scheme in most of
the test cases. For low polynomial degrees and on the coarse
meshes, the difference is often more than a factor of 2, while for
well-resolved solutions, CDG and LDG are similar, and BR2 is about
10\% less accurate. We can also see that all schemes give optimal
convergence rates close to $p+1$ for $||u- u_h||_{0,\Omega}$.

\begin{table}
\begin{center}
\begin{tabular}{ll|rrrrr|c}
\hline
$p$ & Scheme & \multicolumn{1}{c}{$n=2$} & \multicolumn{1}{c}{$n=4$} &
 \multicolumn{1}{c}{$n=8$} & \multicolumn{1}{c}{$n=16$} &
 \multicolumn{1}{c|}{$n=32$} & Rate \\
 \hline
\multicolumn{1}{c}{} & \multicolumn{1}{c|}{} &
\multicolumn{1}{c}{}& \multicolumn{1}{c}{} &
 \multicolumn{1}{c}{} & \multicolumn{1}{c}{}&
 \multicolumn{1}{c}{} & \multicolumn{1}{|c}{}\\[-7pt]
 1 & CDG & $4.54\cdot 10^{-2}$ & $1.52\cdot 10^{-2}$ & $4.62\cdot 10^{-3}$ & $1.25\cdot 10^{-3}$ & $3.27\cdot 10^{-4}$ \ & 1.9 \\
 & LDG & $1.34\cdot 10^{-1}$ & $1.73\cdot 10^{-2}$ & $4.68\cdot 10^{-3}$ & $1.25\cdot 10^{-3}$ & $3.26\cdot 10^{-4}$ \ & 1.9 \\
 & BR2 & $8.60\cdot 10^{-2}$ & $3.08\cdot 10^{-2}$ & $9.23\cdot 10^{-3}$ & $2.47\cdot 10^{-3}$ & $6.36\cdot 10^{-4}$ \ & 2.0 \\
 & & & & & & & \\
 2 & CDG & $8.99\cdot 10^{-3}$ & $1.79\cdot 10^{-3}$ & $2.55\cdot 10^{-4}$ & $3.35\cdot 10^{-5}$ & $4.28\cdot 10^{-6}$ \ & 3.0 \\
 & LDG & $3.81\cdot 10^{-2}$ & $2.92\cdot 10^{-3}$ & $3.03\cdot 10^{-4}$ & $3.59\cdot 10^{-5}$ & $4.42\cdot 10^{-6}$ \ & 3.0 \\
 & BR2 & $1.66\cdot 10^{-2}$ & $2.75\cdot 10^{-3}$ & $3.16\cdot 10^{-4}$ & $3.75\cdot 10^{-5}$ & $4.60\cdot 10^{-6}$ \ & 3.0 \\
 & & & & & & & \\
 3 & CDG & $2.61\cdot 10^{-3}$ & $2.44\cdot 10^{-4}$ & $1.71\cdot 10^{-5}$ & $1.10\cdot 10^{-6}$ & $7.03\cdot 10^{-8}$ \ & 4.0 \\
 & LDG & $5.88\cdot 10^{-3}$ & $3.81\cdot 10^{-4}$ & $2.04\cdot 10^{-5}$ & $1.18\cdot 10^{-6}$ & $7.23\cdot 10^{-8}$ \ & 4.0 \\
 & BR2 & $5.64\cdot 10^{-3}$ & $3.77\cdot 10^{-4}$ & $2.47\cdot 10^{-5}$ & $1.52\cdot 10^{-6}$ & $9.46\cdot 10^{-8}$ \ & 4.0 \\
 & & & & & & & \\
 4 & CDG & $1.09\cdot 10^{-3}$ & $4.52\cdot 10^{-5}$ & $1.56\cdot 10^{-6}$ & $5.14\cdot 10^{-8}$ & $1.63\cdot 10^{-9}$ \ & 5.0 \\
 & LDG & $2.04\cdot 10^{-3}$ & $5.00\cdot 10^{-5}$ & $1.65\cdot 10^{-6}$ & $5.28\cdot 10^{-8}$ & $1.66\cdot 10^{-9}$ \ & 5.0 \\
 & BR2 & $1.30\cdot 10^{-3}$ & $6.22\cdot 10^{-5}$ & $2.05\cdot 10^{-6}$ & $6.57\cdot 10^{-8}$ & $2.07\cdot 10^{-9}$ \ & 5.0 \\
 & & & & & & & \\
 5 & CDG & $3.73\cdot 10^{-4}$ & $9.30\cdot 10^{-6}$ & $1.75\cdot 10^{-7}$ & $2.83\cdot 10^{-9}$ & $4.46\cdot 10^{-11}$ & 6.0 \\
 & LDG & $1.06\cdot 10^{-3}$ & $1.32\cdot 10^{-5}$ & $1.93\cdot 10^{-7}$ & $2.91\cdot 10^{-9}$ & $4.50\cdot 10^{-11}$ & 6.0 \\
 & BR2 & $4.42\cdot 10^{-4}$ & $1.08\cdot 10^{-5}$ & $2.05\cdot 10^{-7}$ & $3.31\cdot 10^{-9}$ & $5.23\cdot 10^{-11}$ & 6.0 \\
\hline
\end{tabular}
\caption{$L_2$ errors in the solution for the model Poisson problem, for
 different polynomial degree, $p$, and mesh size, $n$, using CDG, LDG, and BR2 Schemes.}
\label{tab:L2sol}
\end{center}
\end{table}

\subsection{Spectral radius}\label{sec6.3}

In our next study, we compute the spectral radius
$|\lambda_\mathrm{max}|$ of the discretized matrix and compare the
three methods. The spectral radius of the discretized matrix
determines the magnitude of the timestep if an explicit time
marching solution is sought. In Table~\ref{tab:eig}, we show these
values for each of the simulations in the previous section, scaled
by the factor $(h/p)^2$. Here we have used the consistent switch
with the constant $C_{11}=0$ for the CDG and LDG methods and a
value of $\eta=3$ in the BR2 discretization. We observe that the
CDG and the LDG methods have almost identical spectral radii,
while the BR2 method gives about 50\% larger values. It is
possible that a lower value of the $\eta$ parameter in the BR2
method may reduce the spectral radius. However, in this case
stability may be compromised.

\begin{table}
\begin{center}
\begin{tabular}{ll|rrrrr}
\hline
$p$ & Scheme & \multicolumn{1}{c}{$n=2$} & \multicolumn{1}{c}{$n=4$} &
 \multicolumn{1}{c}{$n=8$} & \multicolumn{1}{c}{$n=16$} &
 \multicolumn{1}{c}{$n=32$} \\ \hline
 1
 & CDG & 153.4 & 157.5 & 159.4 & 159.9 & 160.1 \\
 & LDG & 149.5 & 156.7 & 159.2 & 159.9 & 160.1 \\
 & BR2 & 244.0 & 244.8 & 245.2 & 245.4 & 245.4 \\
 & & & & & & \\
 2
 & CDG & 137.4 & 139.8 & 140.8 & 141.1 & 141.1 \\
 & LDG & 135.1 & 139.5 & 140.7 & 141.1 & 141.1 \\
 & BR2 & 216.1 & 215.5 & 215.3 & 215.1 & 215.1 \\
 & & & & & & \\
 3
 & CDG & 159.9 & 161.3 & 161.8 & 162.0 & 162.0 \\
 & LDG & 159.5 & 161.1 & 161.8 & 162.0 & 162.0 \\
 & BR2 & 244.4 & 244.0 & 243.8 & 243.8 & 243.8 \\
 & & & & & & \\
 4
 & CDG & 198.4 & 200.3 & 201.0 & 201.2 & 201.3 \\
 & LDG & 197.7 & 200.2 & 201.0 & 201.2 & 201.3 \\
 & BR2 & 302.1 & 300.9 & 300.6 & 300.6 & 300.6 \\
 & & & & & & \\
 5
 & CDG & 244.8 & 246.0 & 246.4 & 246.5 & 246.5 \\
 & LDG & 245.1 & 246.0 & 246.4 & 246.5 & 246.5 \\
 & BR2 & 368.5 & 368.4 & 368.4 & 368.4 & 368.4 \\
\hline
\end{tabular}
\caption{The spectral radii of the matrices for the model Poisson
problem, scaled by $(h/p)^2$.}
\label{tab:eig}
\end{center}
\end{table}

\section{Conclusions}
\label{sec7}

We have presented a new scheme for discretizing elliptic operators
in the context of discontinuous Galerkin approximations. The main
advantage of the proposed scheme is its reduced sparsity pattern
when compared to alternative schemes such as the LDG, BR2, or IP
methods. This is important when an implicit solution technique is
required. Compared to the LDG scheme the proposed scheme is
compact, meaning that only degrees of freedom in neighboring
elements are connected. Compared to the BR2 and IP schemes, which
are also compact, the CDG scheme produces a smaller number of
nonzero entries in the off-diagonal blocks and, at the same time,
the nonzero elements in the CDG scheme are amenable to a dense
block matrix storage. Like the alternative approaches, the
proposed scheme converges optimally, and numerical tests indicate
that the accuracy obtained compares well with that of the LDG or
BR2 schemes. An additional potential advantage of the CDG scheme
over the LDG scheme when both schemes are used with minimal
dissipation (i.e.,\ $C_{11} = 0$ in the interior faces) is its
insensitivity to the face ordering.

\bibliography{cdg}
\bibliographystyle{plain}

\end{document}